

\documentstyle{amsppt}
\magnification=1200
\NoBlackBoxes
\def\C{\Bbb C}
\def\P{\Bbb P}

\def\O{\Cal O}
\def\ord{\text{\rm ord}\ \ }
\def\ls{\vskip.25in}
\def\ss{\vskip.15in}
\def\Q{\Bbb Q}

{\tt{991111}}\ss
\centerline{\bf On Semipositivity of Sheaves of Differential Operators}
\centerline{\bf
and the degree of a unipolar $\Q$-Fano variety}
\vskip.25in
\centerline{Z. Ran}
\centerline{Math Dept. UCR, Riverside CA 92521 USA}
\centerline{\tt ziv\@ math.ucr.edu}\ls
\centerline{\bf Abstract}\ss
{
  We consider normal projective n-dimensional varieties $X$ whose anticanonical
divisor class $-K$ is ample and where every Weil divisor is a rational multiple
of $K$. The index $i$ is the largest integer such that $K/i$ exists as a Weil divisor
.
 We show (i) if X has log-terminal singularities and locally,
1-forms on $X_{\text smooth}$ extend holomorphically to a resolution, 
then $(-K)^n \leq
(\max (in,n+1))^n;$ (ii) if the tangent sheaf of X is semistable, then $(-K)^n
\leq (2n)^n.$ The proof is based on some elementary but possibly surprising slope
estimates on sheaves of differential operators on plurianticanonical sheaves.
Unlike previous arguments in the smooth case (Nadel, Campana,
Kollar-Miyaoka-Mori), rational curves and rational connectedness are not used.
Actually, the proof yields the stronger result that bounds as above hold on the
'local degree' of $(X,-K)$, and as such the bound in case (ii) is sharp. 
}
\ls\noindent
{\bf{
$\underline{\text{Contents:}}$\ls\noindent
Introduction\ss\noindent
1. Preliminaries\par
1.1 Differential operators.\par
1.2 The sheaves $D^{i,1}$.\par
1.3 Splitting properties.\par
1.4 Slopes: curve case.\par
1.5 Slopes: higher-dimensional case.\par
1.6 Slopes in an extension.\ss\noindent
2. Slopes of differential-operator sheaves.\par
2.1 First-order estimates for rank 1.\par
2.2 First-order estimates for higher rank\par
2.3 Higher-order estimates.\ss\noindent
3. Positivity of $T_X$.\ss\noindent
4. Conclusion

}}
\ls\ls\noindent
{\bf{Introduction}}\par
Our purpose is to give some boundedness results for $\Q$-Fano
varieties of Picard number 1. We begin with some basic definitions.
Recall that a $\Bbb Q$-Fano variety is by definition a normal
projective variety $X$ such that the anticanonical divisor class
$-K=-K_X$ is $\Bbb Q$-Cartier and ample. For such $X$ we define the
(Weil, resp. Cartier) {\it {index}} $i=i(X)$
(resp. $i_C=i_C(X)$) to be the largest (resp. smallest)
integer such that $K_X/i$ (resp. $i_CK_X$)
exists as a Weil (resp. Cartier) divisor (see [R] for a discussion of Weil
divisors and reflexive sheaves, and also Lemma 2 below)
, the divisor class group
$N(X)$ to be the group of Weil divisors modulo rational equivalence,
and the Picard number $\rho = \rho (X)$ to be the rank of $N(X)$.
When $\rho = 1$, $X$ is said to be {\it {unipolar}}; note that in
this case the singularities of $X$ are automatically $\Q$-factorial.
We will say that $X$ is $t$-factorial, for a natural number $t$,
if for every Weil divisor $L$, $tL$ is Cartier; the smallest such
$t$ (that is, the minimum annihilator of
$N(X)/\text{Pic}(X)\ $ may be called the {\it denominator} of $X$.
Now if $(X,L)$ is any $\Bbb Q$-polarized variety and $p\in X$ a smooth
point, the 'local degree' $$\delta (L,p)$$ 
is defined to be the supremum of real
numbers $\delta$ such that for all rational numbers $d<\delta $
and $k\gg 0$, the natural map
$$H^0(kL)\to kL\otimes \Cal O_X/m_p^{kd},$$ is not injective, 
i.e. such that
$$\{ D\in |kL| : \ord_p(D)\geq kd\}\neq\emptyset $$
and
$$\delta (L)=\inf (\delta (L,p):p\in X \text{smooth}).$$ 
By Riemann-Roch, clearly
$$L^n\leq \delta (L)^n,$$ an inequality which is sometimes strict
(see Example 3.4 below).\par
To state our result we need one more definition. A normal,
Cohen-Macaulay variety $X$ is said to be {\it {1-canonical}} if
for any resolution $\epsilon :X'\to X$,
the differential 
$$\epsilon^*(\Omega_X)\overset d\epsilon \to\rightarrow\Omega_Y$$
factors through a map
$$\epsilon^*(\Omega^{**}_X)\to \Omega_Y.$$
In other words, 1-forms on the smooth part of $X$
lift to holomorphic forms on $Y$. Thus, 1-canonical
is precisely the analogue for 1-forms of the condition
defining canonical singularities.
Note that the 1-canonical condition is automatically satisfied if
$\ \Omega^{**}_X=\Omega_X/$(torsion).
This condition is not too restrictive in view
of the following Lemma, whose proof is given at the end of the
Introduction:
\proclaim {Lemma 0} If $X$ is locally a finite quotient of a complete intersection
$Z$ with $Z$ nonsingular in codimension 2, then $X$ has 1-canonical singularities.
\endproclaim

Our main result is the following 
(we work over $\Bbb C$):
\vskip.2in
\proclaim{Theorem 1}  Let $X$ be an $n$-dimensional unipolar
$\Bbb Q$-Fano variety of index $i$
. Then

(i) if $X$ has log-terminal 1-canonical singularities, we have
$$\delta (-K_X)\leq \max (in,n+1)  ; \tag 0.1 $$

(ii) if the tangent sheaf $T_X=(\Omega_X)^*$ is semistable (with
respect to $-K_X$) (singularities only assumed normal), we have
the sharp estimate
$$\delta (-K_X) \leq 2n. \tag 0.2 $$

\endproclaim

\vskip.2in
\proclaim{Corollary 2} For $X$ log-terminal 1-canonical
$t$-factorial as above , we have
$$(-K_X)^n\leq t^n(n+1)^{2n};$$
if $X$ is locally factorial (e.g. smooth)
 of index $1$, then $$(-K_X)^n\leq (n+1)^n;$$
if $X$ is smooth and not $\Bbb P^n$ or a quadric, then
$$(-K_X)^n\leq ((n-1)(n+1))^n;$$ if $X$ is terminal and $n=3$, then
$$(-K_X)^3\leq 4^6i_C^3.$$
\endproclaim
\demo{proof} Follows from the following Lemma,
also proved at the end of the Introduction (I am indebted to
J. McKernan for pointing our an egregious error in an earlier
statement of this, and for suggesting the correction).
\proclaim{Lemma 3} For $X$ $t$-factorial as above (only assumed normal)
, we have\par
 (i) $i(X)\leq t(n+1)$;\par
 (ii) if $X$
is smooth and $i(X)=n+1$ or $n$ then $X$ is $\Bbb P^n$ or a quadric;\par
(iii) if $n=3$ and $X$ is terminal, then $t\leq i_C$.
\endproclaim

\enddemo

For $X$ a smooth unipolar Fano of index 1, Reid [ReMii] has shown
that $T_X$ is stable. It is unknown whether the same holds without the index 1
hypothesis. Smooth unipolar Fanos with  $i=n-1$ have been classified
by Fujita [Fu]; for those with $i=n-2$ a classification has been
announced by Mukai [Mu]. Interestingly, while the smaller the index, the more difficult
the classification, our bound looks best in small-index cases.
\vskip.1in
\proclaim{Corollary 4} For any $t,n>0$,
the family of $n$-dimensional unipolar $t$-factorial 
Fano varieties $X$ with
log-terminal 1-canonical singularities is bounded.

\endproclaim
\demo{proof} In view of Theorem 1, this follows from a result of Koll\'ar
[Koef] which says that there exists $N=N(k,n)$ such that
if $-kK_X$ is cartier then $-kNK_X$ is very
ample (in the smooth case it also follows from the big Matsusaka
theorem).\qed\enddemo

The existence of a universal bound on the anticanonical degrees of all {\it
smooth} unipolar Fano manifolds of given dimension, and hence their
boundedness as a family, was known previously by works of Campana,
Koll\'ar-Miyaoka-Mori and Nadel, see [K] for an exposition and references.
In particular Koll\'ar-Miyaoka-Mori give the bound $(-K_X)^n\leq (n(n+1))^n$ for $X$
smooth.
The singular case is rather different, and I am grateful to J. McKernan
for a crash course on this. Already in dimension 2 the set of log-terminal
unipolar $\Q$-Fanos, as well as their denominators $t$ and degrees $K^2$, are
unbounded, the simplest example, suggested by McKernan, being the
cone over a rational normal curve; see [KeMcK] for many more examples.
It has been conjectured by Alexeev, and proven by him in dim. 2 [A],
that bounding the log-discrepancy by $\epsilon >0$ yields
a bounded family, and Kawamata [Ka] and Borisov [Bor] have proven that
in dim. 3, assuming $X$ terminal (resp. bounding
the Cartier index) does the same. 
Batyrev [B3] has conjectured that the $\Q-$Fano $n-$folds of
bounded Cartier index form a bounded family.

The bound $(-K_X)^n\leq (n+1)^n$, which would
evidently be sharp (e.g. $X = \P^n)$ has apparently been conjectured for
$X$ smooth; it is false in the non-unipolar case by Batyrev's example [b], namely
$X=\P_{\P^{n-1}} (\Cal O \oplus \Cal O(n-1))$.
\par
\proclaim{Corollary 5} Let $X$ be a smooth Fano $n-$fold with $\rho =1$
which is deformable to one
admitting a K\"ahler-Einstein metric. Then $(-K_X)^n\leq (2n)^n.$\endproclaim

\demo{proof} It is well known that existence of 
a K\"ahler-Einstein metric implies $T_X$ stable, so
we can use Theorem 1 (ii).\enddemo
Apparently, it is generally conjectured that all (smooth) Fano manifolds 
with $\rho =1$ should
admit K\"ahler-Einstein metrics ({\it a fortiori} be deformable
to ones that do). On the other hand, it is easy
to make examples of {\it singular} Fano $n$-folds $X$ with $(-K_X)^n>(2n)^n,$
and we conclude that these all have {\it unstable} tangent sheaves.
\vskip.10in
Our proof resembles others in focusing on the existence of
a section $s \in H^0(-kK)$ having a zero of high order
(roughly $k\delta (-K_X)$) at a general point $p \in X$,
but is more 'elementary' in that rational curves and bend-and-break
are not used. Rather,
the basic idea is to consider
a sheaf $D^j(-kK,\Cal O)$
of differential operators
on a plurianticanonical bundle, $k \gg 0$, showing that this
restricts to a semipositive bundle
on a sufficiently ample and general curve-section $C$ of $X$
(through $p$), provided, roughly, that
$j/k \geq i(n+1)$ (log-terminal 1-canonical case) or
$j/k\geq 2n$ (semistable case). Evaluation on  $s$ reduces the vanishing
order by at most $j$, and this yields a contradiction if $\delta (-K)$ is too large.\par
The paper is organized as follows. Sect. 1 gives various preliminary results,
mostly well known, on differential operators and slopes in general. Sect. 2
develops slope estimates for sheaves of differential operators. In sect. 3
we prove generic positivity of the tangent sheaf $T_X$; the method is
cohomological, and it is here that assumptions on the singularities of
$X$ come into play. The argument is concluded in Sect. 4. For a rough
idea of the proof, the reader may wish
to start with Sect. 4, referring back as necessary.
\vskip.10in
\subheading{Acknowledgement}  I would like to thank Professor Dr. E. Viehweg for pointing
out Lemma 1.1, Professors P. Burchard, J. Koll\'ar and R.K. Lazarsfeld for helful
communications, Professor J. McKernan for several helpful comments,
in particular for pointing out an error in an earlier statement
of Lemma 2 as well as suggesting the correct statement,
and most especially Professor H. Clemens for patiently and
generously sifting through many error-filled versions of this
paper.
In addition, I would like to emphasize my debt to several important papers
by Y. Miyaoka.
\demo{proof of Lemma 0}First, if $Z$ is a local complete intersection with singular
locus $S$ of codimension $\geq 3$, we have an exact conormal sequence
$$0\to r\Cal O_Z\to (n+r)\Cal O_Z\to\Omega_Z\to 0.$$

Now it is standard from this sequence
that $\Omega_Z$ is reflexive: indeed,
denoting by $j:Z-S\to Z$ the inclusion,
we have $\Cal O_Z=j_*j^*\Cal O_Z$ by normality and
$R^1j_*j^*\Cal O_Z=\Cal H^{2}_S(\Cal O_Z)=0$ by depth considerations,
hence $\Omega_Z\to j_*j^*\Omega_Z$ is an isomorphism.
\par
Next, writing locally $X=Z/G$, note that we have, for any $G-$ 
modules $A,B$, split epimorphisms $$A\to A^G, B\to B^G, {\text 
{Hom}} (A,B)^G\to {\text {Hom}} (A^G,B^G),$$ whence a split 
epimorphism $${\text {Hom}} (\Omega_Z,\Cal O_Z)^G\to {\text {Hom}} 
(\Omega ^{G}_Z,\Cal O_Z^G). $$ Applying this twice and using the 
previous result, we get a split epimorphism (where $^+$ 
temporarily denotes Hom$(\cdot,\Cal O_Z^G=\Cal O_X)$) 
$$\Omega_Z^G=(\Omega^{**}_Z)^G\to (\Omega^{G}_Z)^{++},$$ which 
easily implies that $\Omega_Z^G$ is a reflexive $\O_X$-module. 
In view of the natural map $\Omega_X\to \Omega_Z^G$, which is generically 
an isomorphism, it follows that $$\Omega_Z^G=\Omega^{++}_X.$$\par  Now 
consider a  resolution $\epsilon:X'\to 
X,$ and let $Z'$ be a resolution of $X'\times_X Z$, so we have a diagram 
$$\matrix Z' &\to & Z \\ \downarrow 
&&\downarrow \\ X' &\to & X\endmatrix $$ We have a pullback map 
$\Omega_Z^G\subset\Omega_Z\to \Omega_{Z'}. $ On the other hand as $X'$ is 
smooth, we can map $\Omega_{Z'}\to \Omega_{X'}$ by the trace, whence
a map $\Omega^{++}_X\to\Omega_{X'}. $
\qed\enddemo
\demo{proof of Lemma 3}(i)
This argument is  due in the smooth case to Kobayashi-Ochiai, see [K].
We will prove more generally that the largest rational
$r\in\Q$ such that $-K_X=rD$ with $D\ $ Cartier satisfies
$$r\leq n+1.$$
Write
$-K_X=iL$ with $tL$  Cartier. It suffices to prove that
$$i/t\leq n+1.$$ Let
$\epsilon :Y\to X$ be a resolution and write
$$K_Y=\epsilon ^*K_X+E $$ with
$E$ an exceptional $\Bbb Q$-
divisor (not necessarily effective), and $M=\epsilon^*L$ (pullback
As $\Q-$Cartier divisor)  which is
nef and big and $tM$ is integral.On the one hand,
for all integers $0\leq j<n, a>0$
we have by
Kawamata-Viehweg $$H^j(-atM)=0.$$
On the other hand, we have by Serre duality
$$h^n(-atM)=h^0(K_Y+atM).$$
Now take $u$ sufficiently large and divisible and write
$$h^0(u(K_Y+atM))=h^0(u(at-i)M+uE).$$
As $bM$ is effective and nontrivial for $b\gg 0$ and as
there is clearly no nonconstant rational function on $Y$ with poles
only on supp$E$, the latter is 0 whenever $at<i$, hence in this case
$$h^n(-atM)=0.$$
Thus the $n$th degree polynomial $\chi (-atM)$ vanishes for
all integers $a\in (0,i/t)$, hence $i/t\leq n+1$, which proves our
assertion.\par
(ii) See [K].\par
(iii)See [Kea], 6.7.2 (see also 6.11.5 for a more complicated bound
for log-terminal threefolds).\qed\enddemo



\ls
\subheading{1. Preliminaries}
\subheading{1.1  Differential operators}
\vskip.10in
For a normal variety $X$, a torsion-free ${\Cal O}_X$-module $M$ and
a locally free module $N$
we denote by $D^i (M,N)$
or $D^i_X(M,N)$ the sheaf
of (holomorphic) $i$-th order differential operators on $M$ with values in $N$,
i.e.
Hom $(P^i (M),N)$ where $P^i$ denotes the $i$-th principal parts (or jet) sheaf.
For $i=1$ we have an exact sequence
$$\Omega_X\otimes M\to P^1(M)\to M\to 0$$
which is left-injective and locally split over the open set
reg$(M)$ where $M$ is locally
free, and induces
$$0\to \text{Hom}(M,N)\to D^1(M,N)\to \text{Hom}(M,N\otimes T_X)$$
(where $T_X=(\Omega_X)^*)$), which is right-surjective over reg$(M)$,
and in particular induces
an exact sequence (called the Atiyah sequence)
$$0\to \text{Hom}(M,N)\to D^1(M,N)\to G\to 0$$
where $G$ is isomorphic over reg$(M)$ (a fortiori in codimension 1)
to $T_X\otimes M^*\otimes N$
Note that $D^i (M,N)$ forms an ${\Cal O}_X$-bimodule, and in fact
there is a natural map
$$
N \otimes D^i ({\Cal O},{\Cal O}) \otimes M^* \to D^i (M,N)
$$
which is an isomorphism over the open set where $M$ is locally free and $X$ is smooth.
As is well known for $X$ smooth, the action of $T_X$ by
Lie derivative on $K_X$ gives rise to
an identification
$$
D^i ({\Cal O},{\Cal O}) = D^i (K_X, K_X)^{op}, \quad \text{op = opposite bimodule},
$$
hence for $M,N$ locally free
 $D^i (M,N) = D^i (N^*\otimes K_X, M^*\otimes K_X)^{op}$.
\ls
\subheading{1.2 The sheaves $D^{i,1}$}\par
As a convenient intermediate
object for passing from $D^i (M,N)$ to
$D^{i+1}(M,N)$ we will consider the sheaf
$$
D^{i,1} (M, N) = D^1 (P^i (M), N) .
$$
Note the following sequences which are defined for all $M,N$
and exact over reg$(M)$:
$$
\align
0&\to D^{i,1} (M,N) \to D^{i+1,1} (M,N) \to D^1 (S^{i+1}(\Omega_X) \otimes M,N)\to 0 \tag1.1\\
0&\to D^i (M,N) \to D^{1,i} (M,N) \to T_X\otimes D^i (M,N) \to 0 \tag 1.2
\endalign
$$

Combining the evident pairing $D^{i,1} (M,N) \times P^i (M) \to N$ with the
canonical $i$-th order differential operator $M \to P^i (M)$, we get a pairing $D^{i,1} (M,N)
\times M \to N$ which is easily seen to be a differential operator of order $i+1$ on $M$,
hence yields a natural map
$$
D^{i,1} (M,N) \to D^{i+1} (M,N).
$$

By an induction using (1.1) we see easily that
this map is surjective locally
over the open set where $M,N$ are locally free and $X$ is smooth:
indeed over this open set we have an exact diagram
$$\matrix
0&\rightarrow& D^{i,1} (M,N)& \rightarrow& D^{i+1,1} (M,N)&\rightarrow &D^1(S^{i+1}\Omega_X
\otimes M,N)&\rightarrow &0 \\
&&\downarrow &&\downarrow &&\downarrow &&\\
0&\rightarrow& D^{i+1} (M,N)& \rightarrow& D^{i+2} (M,N)&\rightarrow & M^*\otimes N\otimes
S^{i+2}T_X &\rightarrow & 0\endmatrix $$
where the left (resp. right) vertical arrow is surjective by induction
(resp. the Atiyah sequence).\par
 Using the
natural map $H^0 (M) \to H^0 (P^i(M))$, the sequence (1.2) gives rise to a pairing
$$
* : D^{i,1} (M,N) \otimes_{\C} H^0 (M) \to T \otimes N \tag 1.3
$$
which is clearly left $\Cal O_X$-linear. \par

Now as functor in $M$, note that $P^i(M)$ is exact over the smooth
part of $X$, because the
$i$-th neighborhood of the diagonal in $X\times X$ is flat over the smooth
part of $X$. Consequently $D^i(M,N)$ is right-left
exact over the smooth part, i.e. an exact sequence
$$
0\to M' \to M \to M'' \to 0
$$
induces
$$
0\to A'' \to D^i(M,N) \to A' \to 0 \tag  1.4
$$
exact with $A''$ and $A'$ isomorphic in codimension 1 to
$D^i(M'',N)$ and $D^i(M',N)$, respectively.

Note in particular that all of the above exactness statements hold
in a neighborhood
of a generic curve-section.\ls

\subheading{1.3 Splitting properties}\par
As is well known for $M,N$ locally free and $X$ smooth the extension class of
the Atiyah extension $$
0 \to M^* \otimes N \to D^1 (M,N) \to T_X \otimes M^* \otimes N \to 0
$$
as left modules is induced from the Atiyah Chern class of $M$, and in particular
is nontrivial if $M$ has some nontrivial (ordinary) Chern class and $N\neq 0$.
As Viehweg kindly pointed out to me, this remark may be considerably strengthened
if $c_1(M)^n \neq 0$.
First a definition. An
exact sequence of locally free sheaves on $X$:
$$
0\to E \to F \to G \to 0\tag 1.5
$$
is said to be {\it strongly nonsplit}
if
the associated extension element in $H^1(G^*\otimes E)$ does not
lie in the image of $H^1(A)\to H^1(G^*\otimes E)$
for any lower-rank subsheaf $A\subset G^*\otimes E$.
If $E$ is invertible, this means precisely that the extension $F$
does not come from a {\it {locally split}} extension of any lower-rank
subsheaf of $G$.
Notice that if $F$ is of rank 2, i.e. $E,G$ are both invertible,
strongly nonsplit is equivalent to nonsplit.
\par
On the other hand an extension of torsion-free sheaves (1.5) is said
to be {\it quasisplit}
if it lies in the image of the natural map
$$
{\text {Ext}}^1(G',E)\to {\text {Ext}}^1(G,E)$$
for some torsion-free lower-rank quotient $G'$ of $G$;
$G'$ or the corresponding (nontrivial, saturated) subsheaf of $G$ is
called a {\it quasisplitting}. Thus, a non-quasisplit extension is one
that does not come from {\it {any}} extension, locally split or not,
of a lower-rank subsheaf of $G$.
Note that a strongly nonsplit
extension (even of vector bundles) may well be quasisplit, as there
will in general be nonlocally- free subsheaves $G'$ and they will admit
nonlocally split extensions inducing locally split extensions of $G$.\par
 The following result, pointed out by Viehweg, will not be
needed in the sequel, but is good for motivation:

\proclaim{Lemma 1.1}  {\rm (Viehweg)}:  Let $X^n$ be smooth compact, $L$ a line bundle
on $X$ with $c_1(L)^n \neq 0, A \subset \Omega ^1_X$ a
subsheaf of rank $<n$.
Then
the extension class $c_1(L) \in H^1 (\Omega^1_X )$
is outside the image of $H^1(A)$.
Hence the extensions
defining $P^1(L)$
and $D^1(L,\Cal O)$ are strongly nonsplit
\endproclaim

\demo{Proof}
Suppose $c_1(L)$ comes from
an element $\alpha \in H^1(A)$. Then we may represent $c_1(L)$ by a
suitable  $\check{C}$ech cocycle $z$ with values in $A$. As $A$ has rank $<n$, the
cup-power $z^n$, which
represents $c_1(L)^n \in H^n(\Lambda^n\Omega^1_X)$ must vanish
'point-by point'(even as a cocycle), against
our hypothesis $c_1(L)^n\neq 0$.\qed
\enddemo

\ls
\subheading{1.4 \ Slopes: curve case}
\par
We begin by reviewing some definitions
and facts about bundles on curves and their slopes
(see [S],[SB], [MehR] for details).  For a vector bundle
$E$ on a smooth curve, we denote by $\mu(E)$
its slope, i.e.
$$\mu(E)={\text{deg}}(E)/{\text rk}(E),$$
by $\mu^{'}(E)$
the 'shifted' slope
$$\mu^{'}(E)={\text{deg}}(E)/({\text rk}(E)+1) $$
and by $\mu_{\max}(E)$
and $\mu_{\min}(E)$ (resp $\mu^{'}_{\max}(E)$ and $\mu^{'}_{\min}(E)) $ the
largest (resp. smallest) slopes (resp. shifted slopes)
of a subbundle (resp. quotient bundle) of $E$.  As is
well known, the former coincide respectively with the slopes of the first and last associated
gradeds of the Harder-Narasimhan (HN) filtration of $E$.  See {\it op. cit.} for various
basic properties of these invariants. One property we need which is
not mentioned there is behaviour with respect to duality, viz.
$$
\mu_{\min}(E)=\mu_{\max}(E^*).
$$
This can be checked easily.\ls
\subheading{1.5\ Slopes: higher-dimensional case}\par
Now given a (normal projective) variety $X$ we shall
henceforth denote by $C$ a sufficiently general
sufficiently ample curve-section of $X$
(say with respect to a given polarization $H$), and
define slopes of a torsion-free sheaf $E$
on $X$ by restricting on $C$.  The results of
[MeR], which show that an $H-$semistable sheaf of $X$ restricts to
a semistable one on
$C$, imply that these slopes coincide with
those based on $H-$semistabilty
on $X$;
in particular they are independent of the choice of
$C$ and $\mu_{\min}(E)$ coincides with the
slope of the last associated graded of the
$HN$ filtration of both $E$ and $E|_C$., which are compatible
(i.e. the former restricts to the latter).
Also various linear algebra type properties of slopes carry
over from $C$ to $X$,
and two torsion-free sheaves which are isomorphic
in codimension one have the same slopes.
$E$ is said to be {\it generically (semi)positive} if
$\mu_{\min}(E) >(\geq) 0$.
Naturally the slopes of an arbitrary coherent sheaf
are defined to be those of its largest torsion-free quotient.
\ls
\subheading{1.6\ Slopes in an extension}\par
Now we want to give a simple remark concerning the behavior of slopes under
extension.  To this end we introduce the following invariant
$$
a = a(X, H) = \min \{A.H^{n-1}: A \subset X \ \ \text{nontrivial effective
Weil divisor}\}.
$$
Note that for
$X$ unipolar $\Q$-Fano and notations as above) we have
$$a\geq -K_X.H^{n-1}/i\geq -K_X.H^{n-1}/(n+1)$$
as $i\leq n+1$ (Lemma 2).
It will be convenient to abuse notation a bit and assume
$C\sim H^{n-1}$.

\proclaim{Lemma 1.2}  Let $0 \to E \to F \to G \to 0$ be an
extension of
torsion-free sheaves on $X$ with $F$ of rank $r+1$
and $E$ reflexive of rank 1.
Then for any torsion-free quotient $F'$ of $F$ that is not a quasisplitting,
we have

$$
\mu (F') \geq \min ((\frac{a+ E.H^{n-1}}{s+1}+\frac{s}{s+1}\mu_{\min}(G),
s=0,...,r),
\mu_{\min} (G)) .\tag 1.5
$$
\endproclaim

\demo{Proof}
Let $E'\subset F'$ be the saturation of the image of $E$.
Then $G'=F'/E'$ is torsion-free of rank $s$ say, and we have an
exact diagram
$$
\matrix
0 & \to & E &\to & F & \to & G & \to & 0 \\
&&\downarrow&&\downarrow&&\downarrow\\
0 & \to &E'&\to&F'& \to & G' & \to & 0 \\
&&&&\downarrow&&\downarrow \\
&&&&0 && 0
\endmatrix
$$

If $E'=0$ , then
 $G {\overset \sim\to\rightarrow} G'$, and the lower bound
on $\mu (F')$ clearly holds .
Otherwise, $E\to E'$ is injective. Suppose this map is
surjective in codimension 1. Then so is the composite map
to the double dual $E\to E'\to E^{'**}$. But $E^{'**}$ is reflexive,
and an injection of rank-1 reflexive sheaves which is surjective in
codimension 1 is an isomorphism (cf. [R]).
It follows that $E=E'$ which contradicts our
hypothesis that
$F'$ is not a quasisplitting. Thus $E\to E'$ must vanish on some nontrivial effective
divisor and again our estimate follows easily.
\enddemo
\ls
\subheading{2. Slopes of differential-operator sheaves}\ls
The purpose of this section is to give some slope
estimates for sheaves of differential operators, culminating
in Proposition 2.6 below, which is a (generic) semipositivity
result  for operators on suitable plurianticanonical bundles.\ls
\subheading{2.1 First-order estimates for rank 1}\par
The basic idea here is to apply Lemma 1.2 to sheaves of the form
$F=D^1(L,\Cal O)$, where $L$ is an ample line bundle
on our unipolar $\Bbb Q$-Fano  $X$, and $F$ is considered as an
Atiyah extension. The problem is to deal with possible
quasisplittings.

\proclaim{Lemma 2.1}
For any rank-1 torsion-free sheaf $L$ on $X$ with $L.C\neq 0$ we
have $$\mu_{\min}(D^1(L,\Cal O))\geq -L.C+b\tag 2.1 $$ where
$$b=\min (a,\mu^{'}_{\min}(T)).$$
\endproclaim
\demo{proof} To begin with, there exists $s>0$ so that
$L^{[s]}:=(L^{\otimes s})^{**}$ is invertible
$D^1(L,\Cal O)$ is isomorphic in codimension 1 to
(hence has the same slopes as)
$D^1_T(L^{[s]},\Cal O)\otimes L^{[-s+1]}$. Thus there
is no loss of generality in assuming $L$ invertible.

Now suppose given
a torsion-free quotient $F'$ of $D^1(L,\Cal O)$ .
If $F'$ is not  a
quasisplitting, i.e. $L^*$ does not map isomorphically to a
saturated subsheaf of $F'$, then Lemma 2.1 applies
and note that the RHS of (2.1) in this case is just $$-L.C+\min
(a,\mu_{\min} (T))\geq -L.C+b,$$ hence $$
\mu (F')\geq -L.C+b .$$  Otherwise, $F'$
is a quasisplitting, so we get a saturated
subsheaf $T_0\subseteq T$ so that the Atiyah extension splits
over $T_0\otimes L^*$. If $T_0\neq T$, then obviously
$\mu (F')\geq -L.C+b$. If $T_0= T$, the Atiyah extension for
$L$ splits at least over the smooth part of $X$ and in particular
over $C$. This clearly implies that the Atiyah extension for
$L\otimes \Cal O_C$ splits,which implies that $L\otimes \Cal O_C$
is flat, hence $L.C=0$ which is impossible.\qed

\enddemo\ls
\subheading{2.2\ First-order estimates for higher rank}\par
 We  shall next extend
Lemma 2.1 to higher-order
operators.

\proclaim{Lemma 2.2}  If $E$ is a torsion-free generically positive
sheaf on $X$,
then
$$
\mu_{\min} (D^1 (E, \Cal O_X)) \geq \mu_{\min} (E^*) + b
$$
\endproclaim
(Here and elsewhere $D^i$ are considered as left $\Cal O_X$-modules unless otherwise
specified).

\demo{Proof}  If $E$ has rank 1,
the assertion follows directly from
Lemma 2.2.
In the general case rank $E=\rho$ we may use the HN filtration
$E.$ of $E$ and the map
$$
D^1(E_{i-1},\Cal O)\to D^1(E_i,\Cal O)/D^1(E_i/E_{i-1},\Cal O),
$$
(cf. (1.3)), which is an isomorphism in codimension 1,
to reduce to the
case $E$ semistable.
Then use in a similar manner a Seshadri
stable filtration to reduce to the case $E$ stable
(i.e. $E|_C$ stable, cf. [S]). At this point let us fix $C$ and
consider
a suitable unbranched degree- $d$ cyclic cover
$$
\pi : U' \to U
$$
of a tubular neighborhood of $C$ in $X$ where $\pi ^* \det(E)$ admits
a $\rho$-th root $L$,
so we may write
$$\pi ^*E = F \otimes L$$ where $c_1 (F) = 0$. Let $C' = \pi^{-1}(C)$.
Note that $L.C'=dc_1(E).C=d\mu_{min}(E)$. We compute:
$$
\mu_{\min}(D^1_{U'}(L,\Cal O_{C'}))=\mu_{\min}(L^{\rho -1}\otimes D^1
_{U'}(L^{\rho},\Cal O_{C'}))=
$$
$$
=(\rho -1)L.C' + \mu_{\min}(\pi ^*D^1_X(\det E,\Cal O_C))=
(\rho -1)L.C' + d\mu_{\min} (D^1(\det E,\Cal O_X))
\geq
$$
$$
\geq (\rho -1)L.C' + d(\det E^*.C +b) = -L.C' +db
= d(\mu_{\min}(E^*) +b).
$$
(The inequality follows from the rank-1 case)
Now note that $D^1_{U'} (\pi ^*E, \Cal O_{C'})$
depends only on the restriction of
$\pi ^*E$ on the first-order neighborhood $C_1$ of $C'$ in
$U'$, and that $C_1$ coincides with the first-order neighborhood of the
zero-section of the normal bundle $N$ to $C'$ in $U'$, hence admits a
projection map $p:C_1 \to C'$ and a scaling
action, i.e. a family of endomorphisms $\{\phi _t, t\in \C \}$ given
by multiplication by $t$ on a fibre of $N$.
This gives rise to a family of sheaves
$$
\{F_t = \phi _t^* (F) : t \in \C \}
$$
with $F_t \simeq F$ for all $t\neq 0$ and $F_0 =  p^*(F|_{C'})$.
Now by an easy general remark about bundles on curves,
the pullback of a stable bundle by an unramified cyclic cover
is stable. Consequently,
as $E|_C$ is stable, so is $\pi ^*E|_{C'}$, hence also $F|_{C'}$,
so that
$F|_{C'}$
has the form
$\Cal F \otimes \Cal O_{C'}$ for some locally constant sheaf $\Cal F$ on $C'$.
As $\mu_{\min}$ decreases under specialization, we get
the estimate
$$
\mu_{\min} (D^1 (E, \Cal O_C)) =
(1/d) \mu_{\min} (D^1_{U'}(\pi ^*E,\Cal O_{C'})) =
(1/d) \mu_{\min} (D^1_{U'} (F \otimes L, \Cal O_{C'}))
$$ $$
\geq (1/d) \mu_{\min} (D^1_{U'}(F_0 \otimes L, \Cal O_{C'}))
= (1/d) \mu_{\min} (\Cal F^{\surd} \otimes_{\C} D^1
_{U'}(L, \Cal O_{C'}))
$$ $$
= (1/d)\mu_{\min}(D^1_{U'}(L,\Cal O_{C'}))\geq
b+ \mu_{\min} (E^*).
$$
\qed
\enddemo\ls
\subheading{2.3\ Higher-order estimates}\par
Next we extend this slope estimate to higher order,
via the sheaves $D^{m,1}$.  First observe
the standard formula
$$
c_1 (D^m (E, \Cal O_X)) = {\binom{r+m}{r+1}} \rho c_1 (T) +
{\binom{r+m+1}r} c_1 (E^*).
$$
Then a straightforward induction based on Lemma 2.2 plus the surjection $D^{1,m} (E,
\Cal O_X) \to D^{m+1} (E, \Cal O_X)$ yield:

\proclaim{Lemma 2.3}  In the above situation we have
$$
\mu_{\min} (D^{m+1} (E, \Cal O_X)) \geq \mu_{\min} (D^{m,1} (E, \Cal O))
\geq
\min(0,\mu_{\min}
(E^*) + (m+1) b).
$$
\endproclaim
\demo{Proof}
To begin with, the first inequality is immediate from the map
$D^{m,1}(E,\Cal O)\to D^{m+1}(E,\Cal O)$ which is surjective in
codimension 1.
Now by induction on $m$, if $Q$ is any quotient of the
HN filtration $F.$ of
$D^m (E, \Cal O)$, then
$$\mu_{\min} (E^*)+mb \leq \mu (Q) .$$
Now either $\mu (Q)<0$, in which case, $Q^*$
being semistable,
$Q^*$ is generically positive and so by Lemma 2.3
$$\mu_{\min} (D^1 (Q^* ,\Cal O)) \geq \mu_{\min}
(E^*)+(m+1)b; $$
or else $\mu (Q) \geq 0$, in which case clearly $\mu_{\min} (D^1 (Q^* ,\Cal O))
\geq 0$. Using (1.3) as above $F.$ induces a filtration on
$D^{m,1} (E, \Cal O) $ whose quotients are isomorphic in
codimension 1 to the $D^1(Q^*,\Cal O)$, hence satisfy the
above inequality, so the Lemma holds. \qed

\enddemo

Specializing to the case of a line bundle, we conclude
\proclaim{Lemma 2.5}  Let $L$ be a line bundle on $X$ and $m$ any integer satisfying
$$
\gather
m \geq \frac{L.C}{b} - 1\\
\endgather
$$
Then, in the above situation,
 $D^{m,1} (L, \Cal O_X)$ and $D^{m+1} (L, \Cal O_X)$ are generically (left) semipositive
\endproclaim

Specializing further to the case $L = k \det(T)$, we get the following useful estimate

\proclaim{Proposition 2.6}  Let $X$ be a complex unipolar $\Bbb Q$-Fano
variety.  Then
$D^{\alpha k,1} (-kK, \Cal O_X)$ and
$D^{\alpha k+1}(-kK, \Cal O_X)$ are generically
semipositive provided
 $\alpha k$ is an
integer and $$ \alpha \geq \max (-K.C/\mu_{\min}^{'}(T) , -K.C/a) . \tag 2.5 $$
\endproclaim

\demo{proof}
Observe that the $RHS$ of the above inequality on $\alpha$ is just
$\frac{-K.C}{b}$, hence by hypothesis there exists an integer $m$
with $\frac{-K.C}{b} k - 1 \leq m < k (r+1)$. so that the
hypotheses of Lemma 2.5 are satisfied, hence $D^{m,1} (-kK, \Cal
O_X)$ is generically semipositive, which easily implies that so is
$D^{\alpha k, 1} (-kK, \Cal O_X)$ since $\alpha k \geq m$,
hence so is $D^{\alpha k+1}(-kK,\Cal O_X)$.\qed
\enddemo
To take advantage of this result, it is convenient to introduce the
following definition. For a $\Q$-divisor $L$, we define the 'differential
degree'
$$\gamma (L)\in\Bbb R\cup\{\infty\}$$
 as the inf of all $m\in\Q$
such that for all $\alpha>m$ and all $k$ sufficiently large and divisible,
$D^{\alpha k}(kL,\Cal O_X)$ is generically semipositive. Thus
the conclusion of Proposition 2.6
can be rephrased as the estimate
$$\gamma (-K)\leq \max (-K.C/\mu_{\min}^{'}(T) , -K.C/a) . \tag 2.6 $$
\ls
\subheading{3. Positivity of $T_X$}\ls
By definition, a Fano variety $X$ has a tangent sheaf $T_X$
which is positive 'on average'.
The purpose of this section is to show that, with suitable
extra hypotheses, $T_X$ is actually positive on a generic
curve-section.
We will prove the following result, which is apparently
well known in the smooth case
(and which also is the only place where log-terminal 1-canonical
singularities are used).
\proclaim {Proposition 3.} If $X$ is $\Q$-Fano unipolar
with log-terminal
1-canonical singularities then $T_X$ is generically positive.
\endproclaim
\demo{proof} Let $$T_X\to Q\to 0$$ be a quotient of rank $r>0$ and $c_1\leq 0$,
corresponding to a reflexive saturated subsheaf $Q^*\subset\Omega_X^{**}$
 and to a section
$$u\in H^0(Q^{**}\otimes\Omega_X^{**}).$$
Note $r<n$ and set $$M=c_1(Q^*)=(\bigwedge^rQ^*)^{**}.$$
This is a divisorial sheaf which is
either numerically trivial or
ample (i.e. for some $s>0$, $M^{[s]}$ is Cartier and either
numerically trivial or ample; note that $M^{*[s]}=M^{[-s]}$ is then
Cartier as well  ).
We consider resolutions
$$\epsilon : Y\to X$$
with exceptional locus $\bigcup E_i$ a divisor with simple
normal crossings, and with further good properties to be specified.
Write
$$\epsilon
^* (M^{[s]})=:sN_1$$ for a $\Q$-divisor $N_1$
on $Y$. Note that if $M^{[s]}$ is numerically trivial then
so is $N_1$ (and in this case $N_1$ is actually integral, since
numerically trivial divisors are divisible), while if $M^{[s]}$
is ample then $sN_1-\sum e_iE_i$ is ample for $0<e_i\ll 1$,
hence we can write
$$N_1=A_1+F_1$$
where $A_1$
is $\Q$-ample
 and $F_1$ is $\Bbb Q$-
effective with simple normal crossing support and integral part $[F_1]=0$.
Now the following general result is a
sort of 'embedded resolution' for sheaves and is elementary
but worth noting:
\proclaim{Lemma 3.2} Given a torsion-free sheaf $S$ on $Y$, there is
a blowup $\eta :Z\to Y$ such that $\eta^*S/${\rom{(torsion)}} is locally free.
\endproclaim
\demo{proof}
Let $\alpha :F\to S$ be a surjection with $F$ locally free
(in fact of the form $\oplus H^*$ with $H$ sufficiently ample), and
consider the 'canonical resolution' of $\alpha$, i.e.
let $Y'$ be the subvariety of the Grassmann bundle
$$\pi :G(r, F)\to Y,
r ={\text{rank}}(F)-{\text{rank}}(S)$$
 defined as the zero-locus
(with the reduced structure)
of the natural map $Sub\to \pi^*S$,
where $Sub$ denotes the tautological subbundle,
 with natural map
$\beta :Y'\to Y$. Clearly $\beta$ is birational.
Let $Z$ be a desingularization of $Y'$ with natural map $\eta :Z\to Y$,
and $Q$ the pullback of the tautological quotient bundle on $G(r, F)$ to
$Z$. Then there is an induced surjection $Q\to\eta^*S$, and the induced map
$Q\to\eta^*S$/(torsion) is surjective and generically injective, hence
an isomorphism.
\qed\enddemo
It follows in our situation that we may assume
$$-N:=\epsilon^*(M^*)/({\text {torsion}})$$ is
invertible and that $\epsilon^*(Q^{**})/({\text {torsion}})$
is locally free. Note that $$c_1(\epsilon^*(Q^{**})/({\text {torsion}}))=-N.$$
Now the multiplication map
$$(\epsilon ^*(M^*))^{\otimes s}=\epsilon^*((M^*)^{\otimes s})
\to \epsilon^*(M^{[-s]})$$
gives rise to a sheaf inclusion $-sN\subset -sN_1$,
which is an equality locally off the exceptional locus.
Consequently we may write in case $M$ is numerically trivial
$$N=R+S$$
with $R$ numerically trivial and
$S=\sum e_iE_i, e_i\in\Q^{\geq 0}.$
If $M$ is ample we write
 $$N=A+B+F$$
with $A$ $\Q$-ample, $B$ effective integral and $F\ $
$\Q$-effective with normal crossing support and $[F]=0$(i.e. $B+F=S$
as above).

Now by our  hypothesis of 1-canonical singularities the differential of $\epsilon$
factors through a map
$$d\epsilon :\epsilon^*(\Omega_X^{**})\to \Omega_Y$$
and $df(f^*u)$ yields a section $v\in H^0(\Omega_Y\otimes f^*(Q^{**}))$;
by rank considerations the component of Sym$^{2r}v$ in
$H^0(\Omega^r_Y\otimes\bigwedge^rf^*(Q^{**}))$ is nonzero,
whence
a nonzero section
$\O\to \Omega^r_Y(-N)$.


Now our assertion follows from
\proclaim{Lemma 3.3} In the above situation, we have $H^0(\Omega^r_Y(-N))=0$.
\endproclaim
\demo{proof} First, if $M = \Cal O_X$ then
$N=\Cal O_Y$.
By Hodge symmetry it suffices to prove $$H^r(\Cal O_Y)=0.$$
Now by [KaMaMa],Thm. 1-2-5 (log-terminal
Kodaira vanishing), we have $$H^i(\Cal O_X)=0$$ for $i>0,$ and by op. cit.
Thm 1-3-6 (rationality of log-terminal
singularities) we have $$R^i\epsilon_*(\Cal O_Y)=0$$ for
$i>0$. Hence by Leray $H^r(\Cal O_Y)=0$,
as required.

Next, if $M$ is numerically trivial, then as we have just proven
$H^1(\Cal O_Y)=0$,
we have $R=0$, hence $ N$ is $\Q$-effective exceptional.
On the other by definition
$-N$ is effective in a neighborhood of any fibre of $\epsilon$,
and it follows that $N$ is trivial, so we are done as above.

Now assume $M$ is ample.
 Mimicking the usual proof of Kawamata-Viehweg as in [KaMaMa],
Thm. 1-1-1, consider a suitable
finite Galois cover $$\tau :Z\to Y$$ with $Z$ smooth so that
$\tau^*A$ is integral and
$$\Omega^r_Y(-A-B-F)\to\tau_*\tau^*(\Omega^r_Y(-A-B)$$
is a direct summand inclusion, and note the injection
$$\tau^*(\Omega^r_Y)\to \Omega^r_Z.$$ By Nakano,
$$H^0(\Omega^r_Z(-\tau^*A))=0,$$
hence $H^0(\Omega^r_Z(-\tau^*(A+B))=0$,
hence $H^0(\tau^*(\Omega^r_Y(-A-B))=0$,
hence clearly $$H^0(\Omega^r_Y(-A-B-F))=0.\qed$$
\enddemo
\enddemo
\subheading{4.  Conclusion} \vskip .1in We continue with the
notation of Proposition 2.6, and seek firstly to estimate
the RHS of (2.5).
\proclaim{Proposition 4.1} In the above situation we have
(i) if $X$ is log-terminal 1-canonical,
$$\mu ^{'}_{\min}(T)\geq \frac{-K.H^{n-1}}{\max (in,n+1)}; \tag 3.1 $$
(ii) if $T$ is semistable,
$$\mu ^{'}_{\min}(T)\geq \frac{-K.H^{n-1}}{2n}. \tag 3.2 $$
\endproclaim
\demo{proof}
(i) This follows directly from Proposition 3.1 and the definition of index.\par
(ii) Semistability means that any torsion-free quotient $Q$
of $T$ of rank $r>0$ has $\mu (Q)\geq \frac{-K.H^{n-1}}{n}$, hence
$$ \mu ^{'}(Q)=\frac{r}{r+1}\mu (Q)\geq \frac{1}{2}\mu (Q)\geq \frac{-K.H^{n-1}}{2n}.$$
\qed
\enddemo
In light of Propositions 4.1 and 2.6, Theorem 1 follows immediately
from the following easy remark
\proclaim{Lemma 4.2} For any $\Q$-ample divisor $L$, we have
$$\delta (L)\leq\gamma (L).$$
\endproclaim\demo{proof}. If not, pick a rational number
$$\alpha\in (\gamma (L),\delta (L)).$$
Then we
have that for $k$ sufficiently large and divisible,
$$D^{\alpha k} (kL, \Cal O_C)\ \text{ is
semipositive.}$$  Now the evaluation map gives rise to a (left
$\Cal O_C$-linear) map $$ D^{\alpha k} (kL, \Cal O_C) \to H^0
(kL)^*  \otimes \Cal O_C .$$
By definition of $\delta (L)$, we can find, for a suitably small
$\varepsilon > 0$, a nontrivial section $$ s \in H^0 (kL) \quad
\text{with} \quad \text{ord}_p (s) \geq (\alpha + \varepsilon) k ,\
 p\in X \text{ general}.
$$ Projecting $H^0 (kL)^*$ onto $(\C s)^* = \C$, we get a map
$$ \varphi : D^{\alpha k} (kL, \Cal O_C) \to  \Cal
O_C . $$ Choosing $C,s$ sufficiently general mutually, clearly we
can assume $\varphi$ is nontrivial.  Moreover $\varphi$ evidently
factors through $ \Cal O_C ( - \varepsilon kp)$, which
is a negative sheaf, contradicting semipositivity of
$D^{\alpha k} (kL, \Cal O_C)$.  \qed\enddemo
\ls
{\it {Example 4.3}\ \ } Let $X$ be a smooth hypersurface of degree $d\leq n+1$
in $\Bbb P^{n+1}$, $n\geq 3$. Then $X$ is Fano of index $i=n+2-d$ and the tangent
bundle $T_X$ is stable by [PeW]. By Propositions 2.6 and 4.1,
$D^{\alpha k}(-kK_X,\Cal O_C)$ is semipositive for $\alpha\geq 2n$.
On the other hand, given a general point $p\in X$, let $u\in H^0(\Cal O_X(1))$
be the section defining the tangent hyperplane at $p$, and
$s=u^{ik}\in H^0(-kK_X)$, which has order $2ik$ at $p$,showing in
particular that $\delta (-K_X)\geq 2i$, so the estimate (0.2)
is sharp for $d=2$. Further,
choosing things generally enough so that $C$ passes through $p$ but
is not contained in the tangent hyperplane at $p$, we get as above
a nonzero map
$$\varphi :D^{\alpha k}(-kK_X,\Cal O_C)\to \Cal O_C .$$
For any $\alpha < 2i$, this map clearly factors through
$\Cal O_C(-p)$, and consequently $D^{\alpha k}(-kK_X,\Cal O_X)$
is not generically semipositive. In particular, for $d=2$ the estimate
of Propositions 2.6 and 4.1 is sharp as well, i.e. $D^{\alpha k}(-kK_X,\Cal O_X)$
is generically semipositive iff $\alpha \geq 2n$, and thus
$\delta (-K_X)=\gamma (-K_X)$.

Note that this example shows that
it can happen that $L^n<\delta (L)^n.$ This means that
the natural map $H^0(-kK_X)\to -kK_X\otimes
\Cal O_X/m_p^N$ is not necessarily of maximal rank for all $k,N\gg 0$ and
$p\in X$ general.

\vskip1in

\centerline{\bf References}
\vskip.15in
\item{[A]} Alexeev, V.: 'Boundedness of $K^2$ for log surfaces' Internat. J. Math 5,
(1994), 779-810.
\vskip.1in
\item{[B]}  Batyrev, V.V.:  `Boundedness of the degree of multidimensional Fano varieties',
Vestnik MGU (1982), 22-27.
\vskip.10in
\item{[B3]}  Batyrev, V.V.:'The cone of effective divisors of threefolds' Contemp. math.
{\bf {131}}(1989), part 3, 337-352.
\vskip.1in
\item{[Bor]} Borisov, A.: 'Boundedness theorem for Fano log-threefolds' J. Algebraic
Geometry 5, (1996), 119-133.
\vskip.1in
\vskip.10in
\item{[Ka]} Kawamata, Y.: ' Boundedness of $\Q$-Fano threefolds' Contemp. Math
131 (1989), part 3, 439-445. Amer. math. soc.
\vskip.1in
\item{[KaMaMa]} Kawamata,Y., Matsuda, K., Matsuki, K.:' Introduction
to the minimal model program'  Adv. studies in Pure Math. 10 (1987),
283-360.
\vskip.10in
\item{[KeMcK]} Keel, S., McKernan, J.: (to appear)
\item{[Ko]}  Koll\'ar, J.:  'Rational curves on algebraic varieties' Springer 1996.
\vskip.10in
\item{[Koef]}Koll\'ar, J.:  'Effective base point freeness'. Math. Ann. {\bf 296}
595-605 (1994).
\vskip.10in
\item{[Kea]} Koll\'ar, J., et al.: 'Flips and abundance for algebraic threefolds'
As\'erisque 211 (1992).
\vskip.1in
\item{[KoMiMo]} Koll\'ar, J., Miyaoka, Y., Mori, S.:
'Rationally connected varieties'. J. Alg. Geom. 1 (1992), 429-448;
'Rational connectedness and boundedness of Fano manifolds'
. J. Diff. Geom. 36 (1992), 765-769.
\vskip.10in
\item{[MeR]}  Mehta, V., Ramanathan, A.:  `Semi-stable sheaves on projective
varieties and their restriction to curves'.  Math. Ann. {\bf 258} (1982), 213-224.
\vskip.10in
\vskip.10in
\item{[Mi1]}  Miyaoka, Y.: `Deformations of a morphism along a foliation and
applications'.  Proc. Symp. Pure Math {\bf 46} (1987), 245--268.
\vskip.10in
\item{[Mi2]}  $\underline{\hskip.75in}$:  `The Chern classes and Kodaira dimension
of a minimal variety'.  Adv. Studies in Pure Math. {\bf 10} (1981), 449-476.
\vskip.10in
\item{[Mu]} Mukai,S.:'Biregular classification of Fano threefolds
and Fano manifolds of coindex 3', Proc. Nat. Acad. Sci USA 86
(1989), 3000-3002.
\vskip.10in
\item{[PeW]}Peternell, T., Wisniewski,J.:' On stability of tangent bundles
of Fano manifolds with $b_2=1$' (preprint).
\vskip.10in
\item{[ReMi]} Reid, M.:'Canonical 3-folds'. In: A. Beauville, ed.:
'Journ\`ees d'Angers',
North-Holland 1980.
\vskip.10in
\item{[ReMii]}Reid,M.:'Bogomolov's theorem $c_1^2\leq 4c_2$',in: Algebraic Geometry
Kyoto 1977, 623-642, Tokyo, Kinokuniya.
\vskip.10in
\item{[SB]}  Shepherd-Barron, N.I.:  `Miyaoka's theorems...'.  Ast\'erisque
{\bf 211} (1992), 103-114.
\vskip.10in
\item{[S]} Siu, Y.T.: 'Lectures on Hermitian-Einstein metrics for
stable bundles and K\"ahler-Einstein metrics' , DMV 8, Birkhauser 1987.
\vskip.1in



\enddocument